\documentclass[english]{amsart}

\usepackage{babel}
\usepackage{amstext}
\usepackage{amsmath}
\usepackage{amsfonts}
\usepackage{latexsym}
\usepackage{ifthen}
\usepackage{xypic}
\xyoption{all}
\newcommand{\sO}{{\mathcal O}}
\newcommand{\PN}{{\mathbb P}}
\newcommand{\PF}{{\mathbb F}}
\newcommand{\codim}{{\rm codim}}
\newcommand{\Pic}{{\rm Pic}}
\newcommand{\lra}{\longrightarrow}

\newcommand{\KZ}{{\mathbb Z}}
\newcommand{\KQ}{{\mathbb Q}}
\newcommand{\Bs}{{\rm Bs}\,}
\newcommand{\ms}{\xi}
\newcommand{\fib}{{\mathfrak f}}

\newcommand{\tm}{h}
\newcommand{\eps}{\epsilon}

\newcounter{lemma}

\newtheorem{lemma1}[lemma]{\setcounter{equation}{0}}

\newenvironment{lemma}{\begin{lemma1}{\bf Lemma.}}{\end{lemma1}}

\newenvironment{theorem}{\begin{lemma1}{\bf Theorem.}}{\end{lemma1}}
\newenvironment{theorem2}[1]{\begin{lemma1}{\bf Theorem [#1].}}{\end{lemma1}}
\newenvironment{proposition}{\begin{lemma1}{\bf Proposition.}}{\end{lemma1}}

\newenvironment{remark}{\begin{lemma1}{\bf Remark.}\rm}{\end{lemma1}}

\begin{document}

\title{Gorenstein Fano threefolds with base points in the anticanonical system}
\author{Priska Jahnke}
\author{Ivo Radloff}
\address{Mathematisches Institut \\ Universit\"at Bayreuth \\ D--95440 Bayreuth/Germany}
\email{priska.jahnke@uni-bayreuth.de}
\email{ivo.radloff@uni-bayreuth.de}
\date{\today}
\maketitle

%%%%%%%%%%%%%%%%%%%%%%%%%%%%

\section{Introduction}
\setcounter{lemma}{0}

In the classification of Fano varieties, those which are not ``Gino
Fano'', i.e., for which ${-}K_X$ is ample but not very ample, are usually annoying. In the beginning of his classification of Fano threefolds Iskovskikh listed those for which $|{-}K_X|$ is not free. The purpose of this article is to see how his result extends to the canonical Gorenstein case. 

If $X$ is a Gorenstein Fano threefold with at worst canonical singularities, and $\Bs |{-}K_X| \not= \emptyset$, then the rational map defined by $|{-}K_X|$ goes to a surface $W$, which is a rational ruled surface $\Sigma_e$ with $e \ge 0$ or $\widehat{C}_d$, the cone over a rational normal curve of degree $d$. The following Theorem lists the possible pairs $(X, W)$:

\begin{theorem} \label{gen}
Let $X$ be a Gorenstein Fano threefold with at worst canonical singularities and $\Bs |{-}K_X| \not= \emptyset$. Then we are in one of the following cases.
\begin{enumerate}
  \item $\dim \Bs |{-}K_X| = 0$. In this case $X$ is a complete intersection in $\PN(1^4, 2, 3)$ of a quadric $Q$, defined in the first four linear variables, and a sextic $F_6$; $({-}K_X)^3 = 2$ and $W$ is the quadric $Q$ in $\PN_3$.
 \item $\dim \Bs |{-}K_X| = 1$. Then $\Bs |{-}K_X| \simeq \PN_1$ and either 
  \begin{enumerate}
    \item $X$ is the blowup of a sextic in $\PN(1^3,2,3)$ along a complete intersection curve of arithmetic genus $1$; $({-}K_X)^3 = 4$ and $W \simeq \Sigma_1$ or
    \item $X \simeq S_1 \times \PN_1$, where $S_1$ is a del Pezzo surface of degree $1$ with at worst Du Val singularities; $({-}K_X)^3 = 6$ and $W \simeq \PN_1 \times \PN_1$ or
    \item $X = X_{2m-2}$ is an anticanonical model of the blowup of
      the variety $U_m$ (see below) along a smooth, rational complete intersection curve $\Gamma_0 \subset U_{m, \mathrm{reg}}$ for $3 \le m \le 12$; $({-}K_X)^3 = 2m-2$ and $W \simeq \widehat{C}_m$.
    \end{enumerate}
  \end{enumerate}
\end{theorem}
\noindent Here $U_m$ denotes a double cover of $\PN(\sO_{\PN_1}(m)
\oplus \sO_{\PN_1}(m-4) \oplus \sO_{\PN_1})$ with at worst canonical
singularities, such that ${-}K_{U_m}$ is the pullback of the
tautological line bundle $\sO(1)$. For $m \ge 4$, this is a
hyperelliptic Gorenstein almost Fano threefold of degree $4m-8$. The
curve $\Gamma_0$ lies over the complete intersection of some general
element in $|\sO(1)|$ and the ``minimal surface'' $B \in |\sO(1) -mF|$,
where $|F|$ denotes the pencil (note that $\Gamma_0$ is always
contained in the ramification locus). If $m = 3$, then $\Gamma_0$ is
the only curve, on which ${-}K_{U_3}$ is not nef. For details of the construction see section~\ref{examples}.

The cases (a) and (b) are as in Iskovskikh's list. In a different context case i) appears in \cite{Mella} and \cite{IT}, and apparently also in \cite{MoriUtah}.

\vspace{0.2cm}

\noindent {\bf Acknowledgements.} The authors gratefully acknowledge
support by the Schwerpunkt program {\em Globale Methoden in der komplexen Geometrie} of the Deutsche Forschungsgemeinschaft. They also want to thank the referee for many valuable remarks and comments.

%%%%%%%%%

\section{Preliminaries} \label{prel}
\setcounter{lemma}{0}

We recall the following fundamental results:

\begin{theorem2}{Shokurov, \cite{Shokurov}/Reid, \cite{ReidKaw}}
 Let $X$ be a Gorenstein Fano threefold with at worst canonical
 singularities. Then $|{-}K_X|$ contains an irreducible surface $S$
 with at worst Du Val singularities, called {\em general elephant}.
\end{theorem2}

The birational contraction $\tm\colon Y \to X$ in the following theorem is called a {\em partial crepant resolution} or {\em terminal modification} of $X$:

\begin{theorem2}{Reid, \cite{Reid}/Kawamata, \cite{Kaw}}
Let $X$ be a threefold with only canonical singularities. Then there exists a $\KQ$--factorial threefold $Y$ with only terminal singularities and a birational contraction $\tm\colon Y \to X$ such that $K_Y = \tm^*K_X$.
\end{theorem2}
\noindent If $X$ is Gorenstein, then $Y$ is in fact factorial (for
example \cite{Kaw}, Lemma~5.1.). 

A Gorenstein threefold $X$ for which ${-}K_X$ is big and nef is called
{\em almost Fano}. It is called {\em hyperelliptic}, if $|{-}K_X|$ is
free, but the associated map $\varphi$ fails to be injective at
the generic point. In that case  
  \[\varphi\colon X \lra W \subset \PN_N\]
is generically 2--to--1 and $W$ is a so--called variety of minimal degree, i.e.,
  \[\deg W = \codim W + 1.\]
Varieties of minimal degree have been classified by del Pezzo
(\cite{delPezzo}) in dimension $2$ and by Bertini in arbitrary
dimension $n$ (\cite{Bertini}). The list (with some repetitions) is as follows:
\begin{enumerate}
  \item $\PN_n$;
  \item the $n$--dimensional quadric $Q_n \subset \PN_{n+1}$;
  \item (a cone over) the Veronese surface; 
  \item (a cone over) a rational scroll.
\end{enumerate}
The {\em cone over a (rational) scroll}, denoted $\overline{\PF(d_1, \dots, d_n)}$, is the image of 
 \[\PF(d_1, \dots, d_n) = \PN(\sO_{\PN_1}(d_1) \oplus \cdots \oplus \sO_{\PN_1}(d_n)), \quad d_1 \ge  \dots \ge d_n \ge 0\]
in $\PN_{d_1 + \cdots + d_n + n -1}$ under the map associated to the
tautological line bundle which will be denoted $\sO(1)$. Note that for $d_n \ge 1$, $\overline{\PF(d_1, \dots,
  d_n)}$ and $\PF(d_1, \dots, d_n)$ are isomorphic. The pencil on
$\PF(d_1, \dots, d_n)$ will be denoted by $|F|$.

\vspace{0.2cm}

Any effective divisor $D$ on $\PF(d_1, \dots, d_n)$ is in a system
  \[D \in |\sO(k) - lF|, \quad k \ge 0 \mbox{ and } l \in \KZ.\]
Fiberwise, $D \cap F$ is a hypersurface of degree $k$ in $\PN_{n-1}$. If $x_1, \dots, x_n$ denote homogeneous coordinates of $\PN_{n-1}$ corresponding to the summands of our vector bundle, then the monomial $x_1^{e_1} \cdots x_n^{e_n}$ with $e_1 + \cdots +e_n = k$ has as coefficient a function taken from
  \[H^0(\PN_1, \sO_{\PN_1}(e_1d_1 + \cdots + e_nd_n - l)).\]
We will use this in the following form. Consider $\PF(m, m-4) \simeq \Sigma_4$. Denote by $\ms_4$ the minimal section. Any divisor 
  \[D \in |\sO(k) - lF|, \quad  k \ge 0 \mbox{ and } l > k(m-4)\]
contains $\ms_4$ as a component. Indeed, using the above notation, $\ms_4$ corresponds fiberwise to $x_1 = 0$. It therefore suffices to prove that the coefficient function of $x_2^k$ vanishes. This is a section of $\sO_{\PN_1}(k(m-4)-l)$, so the claim follows. 

%%%%%%

\section{The General Elephant in the Case $\Bs |{-}K_X| \not= \emptyset$} \label{elephant}
\setcounter{lemma}{0}

Let $X$ be a canonical Gorenstein Fano threefold with $\Bs |{-}K_X|
\not= \emptyset$. Choose a general elephant $\bar{S} \in |{-}K_X|$. By the
Kawamata--Viehweg vanishing theorem $H^0(X, {-}K_X) \lra H^0(\bar{S},
{-}K_X|_{\bar{S}})$ is surjective, implying
  \[\Bs |{-}K_X| = \Bs |{-}K_X|_{\bar{S}}| \not= \emptyset.\]
Let $\nu\colon S \to \bar{S}$ be a minimal desingularisation of
$\bar{S}$. By Saint--Donat's results on linear systems on smooth K3
surfaces (\cite{SD} or \cite{Shin}),
 \[\nu^*|{-}K_X|_{\bar{S}}| = |\Gamma + mf|,\]
where $m \ge 2$ and
\begin{enumerate}
  \item $|f|$ is an elliptic pencil and
  \item $\Gamma = \Bs |\Gamma + mf| \simeq \PN_1$ is a section.
\end{enumerate}
Let $\Gamma' \subset S$ be an irreducible curve contracted by
$\nu$. Then $(\Gamma + mf)\cdot \Gamma' = 0$, implying $\Gamma \cap
\Gamma' = \emptyset$ or $\Gamma = \Gamma'$. In the first case $S$ and
$\bar{S}$ are isomorphic near $\Gamma$ and $\Bs |{-}K_X| \simeq \PN_1
\subset \bar{S}_{\mathrm{reg}}$. In the second case, $\Gamma$ is
contracted to a point, $\Bs|{-}K_X| = \{p\}$ and $p \in
X_{\mathrm{sing}}$. This is part of a result of Shin:

\begin{theorem2}{Shin, \cite{Shin}} \label{Shin}
Let $X$ be a Gorenstein almost Fano threefold with at worst canonical singularities and assume $\Bs |{-}K_X| \not= \emptyset$. With $\bar{S} \in |{-}K_X|$ a general member we have  
\begin{enumerate}
 \item if $\dim \Bs |{-}K_X| = 1$, then scheme--theoretically $\Bs |{-}K_X| \simeq \PN_1$ is contained in $X_{\mathrm{reg}}$ and $\Bs |{-}K_X| \cap {\rm Sing}(\bar{S}) = \emptyset$;
 \item if $\dim \Bs |{-}K_X| = 0$ then $\Bs |{-}K_X|$ consists of exactly one point and $\bar{S}$ has an ordinary double point at $\Bs |{-}K_X|$. In this case $\Bs |{-}K_X| \subset {\rm Sing}(X)$.
\end{enumerate}
\end{theorem2}
\noindent Note that in the case $\Bs|{-}K_X| = \{p\}$ we have $(\Gamma
+ mf).\Gamma = 0$ on $S$, implying $m = 2$ and hence $({-}K_X)^3 = 2$.

%%%%%%%%%%%%%%%%%%%%

\section{The Case $\dim \Bs |{-}K_X| = 0$} \label{point}
\setcounter{lemma}{0}

Let $X$ be the complete intersection of a quadric $Q$
in the linear variables and a sextic $F_6$ in $\PN(1^4, 2, 3)$. If we
choose $F_6$ general enough, then (see \cite{Mella})
  \[X \cap \{x_0 = x_1 = x_2 = x_3 = 0\} = [0:0:0:0:{-}1:1] = p\]
and $X$ does not meet the singular locus of $\PN(1^4, 2, 3)$. Then $Q$
and $F_6$ are Cartier near $X$ and by adjunction, ${-}K_X \simeq
\sO_{\PN}(1)|_X$ and therefore $\Bs|{-}K_X| = \{p\}$. The rational map defined by $|{-}K_X|$ sends $X$ to the quadric in $\PN_3$ defined by $Q$. 

\begin{proposition}
If $\dim \Bs |{-}K_X| = 0$, then $X$ is as above a complete intersection in $\PN(1^4, 2, 3)$ of a quadric $Q$, defined in the first four linear variables, and a sextic $F_6$.
\end{proposition}

\begin{proof}
(See \cite{Mori}, \cite{Mella}, \cite{IT}). We know $({-}K_X)^3 = 2$
(see the last section). By the Riemann--Roch theorem we get $h^0({-}K_X) = 4$. Let 
  \[x_0, \dots, x_3 \in H^0({-}K_X)\]
be generating sections. We have $h^0({-}2K_X) = 10 = \dim S^2H^0({-}K_X)$. But $|{-}2K_X|$ is base point free, so there exists some 
  \[y \in H^0({-}2K_X), \quad y \not\in S^2H^0({-}K_X).\]
Then we must have a nontrivial relation $Q$ in $S^2H^0({-}K_X)$. The
$x_i$ and $y$ then define a $20$--dimensional subspace of
$H^0({-}3K_X)$. By the theorem of Riemann--Roch $h^0({-}3K_X) =
21$. Denote the missing function by $z \in H^0({-}3K_X)$. Continuing
in this way, we see that there must be a nontrivial relation $F_6$ in
$H^0({-}6K_X)$. In the end $X$ is the complete intersection of $Q$ and $F_6$ in $\PN(1^4, 2, 3)$.
\end{proof}

\begin{remark} \label{remcone}
Since $Q$ is singular at $p$, any $S \in |{-}K_X|$ is singular at $p$. If we choose $Q$ and $F_6$ general, $p$ will be a terminal point of $X$. If we take for $Q$ the quadric cone, $X$ will have canonical singularities along a curve.
\end{remark}

%%%%%%%%%%%%%%%%

\section{The Examples for the Case  $\dim \Bs |{-}K_X| = 1$} \label{examples} 
\setcounter{lemma}{0}

Let $U$ be a canonical Gorenstein threefold. Assume that $|{-}K_U|$
contains a smooth K3 surface $S$ such that
  \[-K_U|_S = 2\Gamma_0 + mf\]
for some $m \ge 3$. Here $\PN_1 \simeq \Gamma_0 \subset
U_{\mathrm{reg}}$ and $|f|$ is an elliptic pencil as in
section~\ref{elephant}. Note that $U$ is a hyperelliptic almost Fano threefold for $m \ge 4$.

Let $Y = Bl_{\Gamma_0}(U)$ be the blowup of $U$ in $\Gamma_0$. The strict transform of $S$
is a smooth K3 surface in $|{-}K_Y|$ which we denote by $S$ as well. We have
  \[-K_Y|_S = \Gamma_0 + mf,\]
implying $\Bs|{-}K_Y| = \Gamma_0 \simeq \PN_1$. An anticanonical model $X$ of
$Y$ is a canonical Gorenstein Fano threefold for which $\Bs|{-}K_X| \simeq \PN_1$.

Examples for $U$ as above are constructed as follows. For $m \ge 4$, $U$ is
almost Fano and the anticanonical map associated to $-K_U$ sends $U$
to a variety of minimal degree
  \[U \lra W \subset \PN_{2m-2}.\]
Here $S$ is sent to $\Sigma_4$, the fourth Hirzebruch surface. The idea
is therefore to construct $U$ as a ramified twosheeted covering of
some variety of minimal degree, for which a general hyperplane section is
isomorphic to $\Sigma_4$. 

We now come to the examples in ii) in reverse order.

\vspace{0.2cm}

\noindent {\bf Examples ii), (c).} The projective bundle
 \[W = \PF(m, m-4, 0), \quad m \ge 3\]
is a resolution of a cone over $\Sigma_4$. The projection of the underlying bundle onto the first two summands gives a split exact sequence and a smooth surface in $|\sO_W(1)|$ isomorphic to $\Sigma_4$. For simplicity, we denote it by
  \[\Sigma_4 \in |\sO_W(1)|.\]
There exists a unique section $B \in |\sO_W(1) - mF|$ meeting $\Sigma_4$ in its minimal section $\ms_4$. Below we prove that for $m \le 12$ we may choose 
 \[D \in |\sO_W(4) - (4m-12)F|,\]
such that the square root of $D$ yields a threefold $U_m$ with at
worst canonical singularities. We have 
 \[\mu\colon U_m \stackrel{2:1}{\lra} \PF(m, m-4, 0) \quad \mbox{and} \quad {-}K_{U_m} = \mu^*\sO_W(1).\]
The section $\ms_4 = B \cap \Sigma_4 \subset D_{\mathrm{reg}}$. Its
reduced inverse image in $U_m$ will be denoted by $\Gamma_0$. As in
ii) (c) of Theorem~\ref{gen}, we denote by $X_{2m-2}$ an anticanonical
model of $Bl_{\Gamma_0}(U_m)$ for $3 \le m \le 12$. We claim that
$X_{2m-2}$ are canonical Gorenstein Fano threefolds with
base locus $\Bs|{-}K_{X_{2m-2}}| \simeq \PN_1$.

\vspace{0.2cm}

In order to prove this it suffices to show that for $D$ general enough
each $U_m$ is a canonical Gorenstein threefold as in the beginning of this
section. Since $\Sigma_4$ comes from a splitting sequence, $D \cap \Sigma_4$ is a general member of
  \[|4\ms_4 + 12\fib|,\]
with $\fib \simeq \PN_1$ a fiber of $\Sigma_4$. A general member of
$|4\ms_4 + 12\fib|$ splits as $\ms_4 + C$ with $C \in |3\ms_4 + 12\fib|$
smooth and disjoint from $\ms_4$ (cf. section~\ref{prel}). The double covering of $\Sigma_4$ yields a smooth K3 surface
$S \in |{-}K_{U_m}| = |\mu^*\sO_W(1)|$ with
  \[\mu_S\colon S \lra \Sigma_4\]
ramified along $\ms_4$ and $C$. The pullback of $\fib$ gives an
elliptic pencil $|f|$ on $S$ with the section $\Gamma_0$ lying over
$\ms_4$ and ${-}K_{U_m}|_S = \mu_S^*\sO(1) = 2\Gamma_0 + mf$. It
remains to show that $U_m$ has at worst canonical singularities for $3 \le m \le 12$ and $\Gamma_0 \subset U_{m, \mathrm{reg}}$.

\vspace{0.2cm}

For $m = 3$ we can choose $D$ and hence $U_m$ smooth and there is nothing to prove. For $m \ge 4$, we always have 
  \[D = B + R\]
with $R \in |\sO_W(3) - (3m-12)F|$. Fiberwise $D \cap F$ consists of a line together with some cubic.

For $4 \le m \le 12$ we can take $R$ to be irreducible, i.e., $D \cap F$ consists of a line and an irreducible cubic. For $m = 4$, the cubic is smooth, meeting the line transversally in
three points. For $m \ge 5$, the line and the cubic intersect in one
point, i.e., in a flex if the cubic is smooth. This gives an A--D--E singularity in the fiber, implying that
$U_m$ indeed has at worst canonical singularities for $3 \le m \le
12$. Since $R.\ms_4 = 0$ we can choose $R$ disjoint from $\ms_4$. Hence $\Gamma_0 \subset U_{m,\mathrm{reg}}$.

For $m \ge 13$ on the other hand, $R = R_1 + R_2 + R_3$ with $R_i \in |\sO_W(1) - (m-4)F|$, so $D \cap F$ consists of four lines through a point. This means that over $F$ we will not have Du Val singularities, implying that $U_m$ is not canonical for $m \ge 13$.

\begin{remark}
The construction works for $m = 2$ as well. Here $\Bs|{-}K_{X_2}| =
\{p\}$ and we get a special case of the threefold $X$ in
section~\ref{point} with $Q$ the quadric cone (see Remark~\ref{remcone}).
\end{remark}

\vspace{0.2cm}

\noindent {\bf Example ii), (b).} The product of $S_1$, a del Pezzo surface with canonical singularities of degree $1$, and $\PN_1$ is a classical example (\cite{Isk}). Choose $8$ points on $\PN_2$ general enough, such that the blowup $\hat{\PN}_2$ of $\PN_2$ in these points still has a nef anticanonical system, and denote by $S_1$ an anticanonical model of $\hat{\PN}_2$. Then $|{-}K_{S_1}|$ is one dimensional by the Riemann--Roch theorem, its members corresponding to elliptic curves passing through the eight points. These curves will meet in a ninth point, implying
  \[\Bs|{-}K_{S_1}| = \{p\}.\]
Then the product $X = S_1 \times \PN_1$ is a canonical Gorenstein Fano threefold with $\Bs|{-}K_X| \simeq \PN_1$.

\vspace{0.2cm}

\noindent {\bf Example ii), (a).} The blowup $X$ in the intersection of two members of $|{-}\frac{1}{2}K_U|$ of the double cover $U$ of the Veronese cone $W$, ramified along a cubic, is a classical example (\cite{Isk}). We give some details to show the connection to the above description.

The blowup of the Veronese cone in its vertex $O$ yields
 \[\PN(\sO_{\PN_2} \oplus \sO_{\PN_2}(2)) \lra W.\]
The strict transform of a special hyperplane section through $O$ gives a $\PN_1$--bundle over a conic. It either decomposes into two copies of $\Sigma_2$ or gives one irreducible surface $\Sigma_4$.

The image of $\Sigma_4$ in $W$ gives $\widehat{C}_4$, the cone over the rational normal curve of degree $4$. In $U$, lying over $\widehat{C}_4$ we find a singular K3 surface $\bar{S} \in |{-}K_U|$ with a double point over $O$. In the reducible case, the two copies of $\Sigma_2$ induce $H_i \in |{-}\frac{1}{2}K_U|$ for $i=1,2$, and their intersection with $\bar{S}$ is the singular point. 

In the blowup $X$ of $U$ along $H_1 \cap H_2$ the singularity of $\bar{S}$ is resolved, i.e., we get a smooth K3 surface $S \in |{-}K_X|$. The same formulas as above show
  \[{-}K_X|_S = \Gamma + 2f\]
with $\Gamma$ the ${-}2$--curve over the singularity and $|f|$ the
induced elliptic pencil. If we choose $H_1, H_2$ general enough, then
$X$ will be a canonical Gorenstein Fano threefold with $\Bs|{-}K_X|
\simeq \Gamma \simeq \PN_1$.

%%%%

\section{The General Setting in the Case $\dim \Bs |{-}K_X| = 1$} \label{gensetting}
\setcounter{lemma}{0}

(cf. \cite{Isk}, \cite{AG5}) By Shin's Theorem, $\Gamma = \Bs |{-}K_X|
\simeq \PN_1 \subset X_{\mathrm{reg}}$. We can write
 \begin{equation} \label{ab}
   N_{\Gamma/X} = \sO_{\PN_1}(a) \oplus \sO_{\PN_1}(b), \quad a \ge b,
 \end{equation}
for some $a, b \in \KZ$. A general elephant $\bar{S} \in |{-}K_X|$ may
have double points, but $\Gamma \subset \bar{S}_{\mathrm{reg}}$. If
$\nu\colon S \to \bar{S}$ denotes a resolution of the singular locus,
then $\nu^*({-}K_X) = \Gamma + mf$, $m \ge 3$, with $|f|$ an elliptic pencil and $\Gamma$ a section
(section~\ref{elephant}). The numbers are related as follows: 
 \[{-}K_X\cdot \Gamma = m-2 = a+b+2.\] 

Let $\sigma\colon X_{\Gamma} \to X$ be the blowup of $X$ along
$\Gamma$ with exceptional divisor $E_{\Gamma} = \PN(N^*_{\Gamma/X}) =
\Sigma_{a-b}$. Then $|{-}K_{X_{\Gamma}}| = |\sigma^*({-}K_X) -
E_{\Gamma}|$ is free, defining a map onto some surface $W$ (\cite{ReidKaw}):
 \begin{equation} \label{xgamma}
   \xymatrix{X_{\Gamma}\ar[d]_{\sigma} \ar[r]^{\varphi} & W & \hspace{-1cm} \subset \PN_{m+1} \\
  X \ar@{-->}[ru] && }
 \end{equation}
The surface $W$ is of minimal degree, i.e., 
  \[m = \deg(W) = \codim(W) + 1.\]
Again by del Pezzo's theorem, in our situation $W$ is one of the following:
\begin{enumerate}
   \item $\widehat{C}_m$, the cone over a rational normal curve of degree $m = a+b+4 \ge 2$, 
   \item $\Sigma_{a-b}$, $a \ge b$.
\end{enumerate}
The map $E_{\Gamma} \to W$ is either an isomorphism or the contraction
of the minimal section. The map $X_{\Gamma} \to W$ is (generically) an
elliptic fibration, and since ${-}K_X$ is ample, any fiber over a
point in $W_{\mathrm{reg}}$ is an irreducible, generically reduced
curve of arithmetic genus one. We distinguish two cases.

\vspace{0.2cm}

\noindent {\bf The case $W$ a smooth ruled surface.} Here we denote by $F_{\Gamma}$ the pullback to $X_{\Gamma}$ of a fiber of $W$, and by $Z_{\Gamma,X}$ the pullback of the minimal section (or the second ruling in the case $W = \PN_1 \times \PN_1$). Note that $|F_{\Gamma}|$ descends to a pencil $|F|$ on $X$. Adjunction on $E_{\Gamma}$ shows
  \[{-}K_{X_{\Gamma}} = Z_{\Gamma,X} + (a+2)F_{\Gamma}.\]
Since $\Gamma \subset X_{\mathrm{reg}}$ and $Z_{\Gamma,X}$ meets $E_{\Gamma}$ transversally near the minimal section $\ms_{a-b}$ of $E_{\Gamma}$, $Z_{\Gamma,X}$ is smooth near $Z_{\Gamma,X} \cap E_{\Gamma}$, and $\sigma(Z_{\Gamma,X}) \simeq Z_{\Gamma,X}$ is smooth near $\Gamma$.

\vspace{0.2cm}

\noindent {\bf The case $W$ a cone.} Here we denote by $F_{\Gamma}$
the strict transform in $X_{\Gamma}$ of a line in $W$ through the
vertex $O$. Notice that this is just a Weil divisor. Let
  \begin{equation} \label{Q}
    \tm'\colon X'_{\Gamma} \lra X_{\Gamma}
  \end{equation}
be a $\KQ$--factorialization of $X_{\Gamma}$ with respect to
$F_{\Gamma}$ (\cite{Kaw}). The map $\tm'$ is small, $X'_{\Gamma}$ is
again Gorenstein with at worst canonical singularities, and the strict
transform $F'_{\Gamma}$ of $F_{\Gamma}$ is $\KQ$--Cartier. We can
choose $X'_{\Gamma}$ such that $F'_{\Gamma}$ is $\tm'$--ample
(\cite{Kaw}). Since $\Gamma \subset X_{\mathrm{reg}}$, both
$X'_{\Gamma}$ and $X_{\Gamma}$ are isomorphic near $E_{\Gamma}$. We
denote the pullback of $E_{\Gamma}$ to $X'_{\Gamma}$ by
$E'_{\Gamma}$. We claim (cf. \cite{Cheltsov})

\begin{lemma}
  On $X'_{\Gamma}$, two general members of $|F'_{\Gamma}|$ do not intersect.
\end{lemma}

\begin{proof}
 Assume $F'_{\Gamma, 1} \cap F'_{\Gamma, 2} \not= \emptyset$. The
 intersection clearly is in the fiber over the vertex $O$ of $W$. Choose an irreducible curve $C \subset F'_{\Gamma, 1} \cap F'_{\Gamma, 2}$. On the one hand, the restriction of some multiple of $F'_{\Gamma, 2}$, which is Cartier, gives an effective Cartier divisor on $F'_{\Gamma, 1}$ supported in the fiber over $O$, implying
  \[F'_{\Gamma, 2} \cdot C \le 0.\]
On the other hand, since $F'_{\Gamma, 1}$ and $F'_{\Gamma, 2}$ do not meet on $E'_{\Gamma}$, we have $C \cap E'_{\Gamma} = \emptyset$. Since ${-}K_{X'_{\Gamma}}\cdot C = 0$ and $E'_{\Gamma}\cdot C = 0$ imply ${\tm'}^*\sigma^*({-}K_X) \cdot C = 0$, the curve $C$ must be $\tm'$--exceptional. Then, by our choice of $X'_{\Gamma}$,
  \[F'_{\Gamma, 2} \cdot C > 0.\]
Hence $F'_{\Gamma, 1} \cap F'_{\Gamma, 2} = \emptyset$.
\end{proof}

\vspace{0.2cm}

Denote by $Y_{\Gamma}$ a terminal modification of $X'_{\Gamma}$. The
pullback of $F'_{\Gamma}$ to $Y_{\Gamma}$ defines a pencil on
$Y_{\Gamma}$, showing that the map to $W$ factors over the blowup
$\Sigma_{a-b}$ of $W$ in $O$. Near $E'_{\Gamma}$, $Y_{\Gamma}$ and
$X'_{\Gamma}$ are isomorphic, and we can blow the divisor down to
obtain $Y$, a terminal modification $\tm\colon Y \to X$ of $X$. We call the map $Y_{\Gamma} \to Y$ again $\sigma$ and end up with the following diagram:
 \begin{equation} \label{diag}
   \xymatrix{ & & \Sigma_{a-b} \ar[d]\\
  Y_{\Gamma} \ar[r] \ar[d]^{\sigma} \ar[urr]^{\psi} & X_{\Gamma} \ar[d] \ar[r] & W\\
  Y \ar[r]^{\tm} & X}
 \end{equation}

Below, we will study $Y$ instead of $X$ and think of $X$ as an
anticanonical model. Note that we have chosen $Y$ as a terminal
modification of a particular $\KQ$--factoriali\-zation of $X$.

For simplicity, denote divisors on $Y_{\Gamma}$ and $X_{\Gamma}$ by
the same letters: the exceptional divisor of $Y_{\Gamma} \to Y$ is
again $E_{\Gamma}$, the curve $\Bs |{-}K_Y| = \Gamma$. The pullback of
a general fiber of $\Sigma_{a-b}$ to $Y_{\Gamma}$ is $F_{\Gamma}$. By
$Z_{\Gamma} + B_{\Gamma}$ we denote the pullback of the minimal
section of $\Sigma_{a-b}$ to $Y_{\Gamma}$, where $Z_{\Gamma}$ denotes
here the unique irreducible component that meets $E_{\Gamma}$ in its
minimal section, and $B_{\Gamma}$ consists of the remaining components,
disjoint from $E_{\Gamma}$. As above we get
 \begin{equation} \label{zbf}
   {-}K_{Y_{\Gamma}} = Z_{\Gamma} + B_{\Gamma} + (a+2)F_{\Gamma}.
 \end{equation}
The pencil $|F_{\Gamma}|$ again descends to the pencil $|F|$ on $Y$. The surface $Z_{\Gamma}$ is smooth near $E_{\Gamma} \cap Z_{\Gamma}$; we will denote the isomorphic images of $Z_{\Gamma}$ and $B_{\Gamma}$ in $Y$ by $Z$ and $B$.

\begin{remark} \label{Kod}
 The general member of the pencil $|F_{\Gamma}|$ is a smooth surface with a relatively minimal elliptic pencil. The intersection $F_{\Gamma} \cap (Z_{\Gamma} + B_{\Gamma})$ is hence either smooth or one of Kodaira's exceptional fibers.
\end{remark}

%%%%%%%%

\section{The Case $W$ a Cone}
\setcounter{lemma}{0}

\begin{proposition} \label{cone}
If $W$ is a cone, then $3 \le m \le 12$ and $X = X_{2m-2}$ is one of
the threefolds constructed in section~\ref{examples}, Examples ii), (c). Here $W = \widehat{C}_m$.
\end{proposition}

\begin{proof}
We use the notation from the last section. Since ${-}K_{X_{\Gamma}}$ is not ample on $E_{\Gamma}$, $b = -2$ and $a \ge 1$ in (\ref{ab}). We can hence use $a+b = m-4$ to eliminate $a$ and $b$ and write everything in terms of $m$:
  \[N_{\Gamma/X} = \sO_{\PN_1}(m-2) \oplus \sO_{\PN_1}({-}2), \quad m \ge 3,\]
and $W = \widehat{C}_m$. In diagram (\ref{diag}), the map from $Y_{\Gamma}$ to $\widehat{C}_m$ now factors over $\Sigma_m$. 

\vspace{0.2cm}

We first assume that $Z$ is $\tm$--nef and show that in this case $Y$
is obtained by blowing up some Gorenstein threefold $V$ along some
smooth curve $\Gamma_0 \simeq \PN_1 \subset V_{\mathrm{reg}}$, such
that $Z$ is the exceptional divisor. We compute
  \[Z \cdot\Gamma = -2 \quad \mbox{and} \quad -K_Y\cdot\Gamma = m-2 > 0.\]
Hence $[\Gamma]$ is contained in the $K_Y$--negative part of $\overline{NE}(Y)$. This part is polyhedral, spanned by $K_Y$--negative extremal rays. The divisor $Z$ is negative on $[\Gamma]$ and nonnegative on any $K_Y$--trivial curve by assumption. We conclude that $Z$ must be negative on at least one extremal ray. Let
 \begin{equation} \label{phi}
   \phi\colon Y \lra V
 \end{equation}
be the contraction of this ray. By \cite{Ben}, the contraction is divisorial, contracting $Z$ either to a curve or to a point. We claim

\begin{lemma} \label{V}
The map $\phi\colon Y \to V$ in (\ref{phi}) is the blowup of a smooth rational curve $\Gamma_0 \subset V_{\mathrm{reg}}$ with normal bundle
  \[N_{\Gamma_0/V} = \sO_{\PN_1}({-}2) \oplus \sO_{\PN_1}(m-4).\]
The contraction is in direction of $|F|$. There exists a smooth K3 surface $S \in |{-}K_V|$, such that
    \[{-}K_V|_S = 2\Gamma_0 + mf\]
with $|f|$ an elliptic pencil induced by $|F|$, and $\Gamma_0 \simeq \PN_1$ a smooth section.
\end{lemma}

\begin{remark} \label{hyp}
The threefold $V$ is a hyperelliptic Gorenstein almost Fano threefold of degree $({-}K_V)^3 = 4m-8$ for $m \ge 4$. For $m = 3$, the anticanonical system is nef on any curve $\not= \Gamma_0$, while
  \[{-}K_V\cdot \Gamma_0 = m-4 = -1.\]
For the case $ m=3$ (as well as $m = 2$) see also \cite{DPS}.
\end{remark}

\begin{proof}[Proof of Lemma~\ref{V} and Remark~\ref{hyp}]
Since $Z_{\Gamma}$ meets $E_{\Gamma}$ transversally in the minimal section, we have $\Gamma \subset Z_{\mathrm{reg}}$. We compute
 \begin{equation} \label{nb}
   \deg N_{\Gamma/Z} = Z_{\Gamma} \cdot_{Y_{\Gamma}}E_{\Gamma}^2 = m-2 > 0.
 \end{equation}
Let us first show $Z \not\simeq \PN_1 \times \PN_1$. If $Z \simeq
\PN_1 \times \PN_1$, then $B \not= 0$, implying that $B$ meets $Z$ in some curve. By (\ref{nb}) $\Gamma$ is ample on $Z$. Then $\Gamma \cap B \not= \emptyset$, which is impossible since $B$ maps to $X_{\mathrm{sing}}$, while $\Gamma \subset X_{\mathrm{reg}}$.

\vspace{0.2cm}

If $Z$ is mapped to a point, then by \cite{Cu}, $Z \simeq \PN_2$,
$\PN_1 \times \PN_1$ or the quadric cone. Since $Z$ comes with a
pencil and $Z \not\simeq \PN_1 \times \PN_1$, all these cases are
impossible. By \cite{Cu}, $Y = Bl_{\Gamma_0}(V)$ the blowup of $V$ in
some curve $\Gamma_0 \subset V_{\mathrm{reg}}$, which is locally a
complete intersection. From $\deg N_{\Gamma/Z} = m-2 > 0$ we conclude that $\Gamma$ maps surjectively onto $\Gamma_0$, and from $\Gamma \subset Z_{\mathrm{reg}}$ we infer that $\Gamma_0$ must be smooth. Then
  \[Z = \PN(N_{\Gamma_0/V}^*) \simeq \Sigma_e \quad \mbox{for some $e > 0$,}\]
where $e > 0$ follows from $Z \not\simeq \PN_1 \times \PN_1$. It is
now clear that $\phi$ is in direction of $|F|$, i.e., fiberwise $\phi$
contracts a ${-}1$--curve in $F$. Denote the induced pencil on $V$
by $|F_V|$. Notice that $Z \simeq \Sigma_e$ implies $B \not= 0$.

\vspace{0.2cm}

1.)~Any curve in $Z_{\Gamma} \cap B_{\Gamma}$ is contracted by
$Y_{\Gamma} \to X_{\Gamma}$, and therefore $B$ intersects $Z$ set
theoretically in the minimal section $\ms_e$ of $Z = \Sigma_e$. Since
$\Gamma$ does not meet $\ms_e$, we conclude $\Gamma = \ms_e +
(m-2)\fib_e$, where $\fib_e$ is a fiber of $\Sigma_e$. From $\Gamma
\cdot_{Z}\Gamma = m-2$ (\ref{nb}) we infer $e = m-2$. Moreover,
$-K_Y\cdot \ms_e = 0$ implies
  \[N_{\Gamma_0/V} = \sO_{\PN_1}({-}2) \oplus \sO_{\PN_1}(m-4).\]
By the adjunction formula, ${-}K_V\cdot\Gamma_0 = m-4$, hence $({-}K_V)^3 = 4m-8$.

\vspace{0.2cm}

2.)~Let $S \in |{-}K_Y|$ be general. Since $S$ meets $Z$ transversally
in $\Gamma$, its image in $V$ is a special member of
$|{-}K_V|$. Identifying $S$ with its image in $V$ we find
  \[{-}K_V|_S = 2\Gamma_0 + mf,\]
where $|f|$ is an elliptic pencil and $\Gamma_0$ is a section (see section~\ref{examples}). If $C \subset V$ is an irreducible curve such that ${-}K_V\cdot C < 0$, then $S\cdot C < 0$ and $C \subset S$. Then ${-}K_V\cdot C = (2\Gamma_0+mf)\cdot C < 0$ so that $\Gamma_0\cdot C < 0$ and hence $C = \Gamma_0$, $m = 3$.
\end{proof}

\vspace{0.2cm}

The argument before Lemma~\ref{V} showing the contractibility of $Z$ in $Y$ requires
$Z$ being $\tm$--nef. In order to achieve this we might have to change
the terminal modification by running the relative $(K_Y + \eps Z)$--program, $\eps \in \KQ^+$, $\eps \ll 1$, with respect to $\tm\colon Y \to X$.

The contraction of any $(K_Y + \eps Z)$--negative extremal ray in $\overline{NE}(Y/X)$ is small; the curves contracted are $K_Y$--trivial and contained in $Z$. After finitely many flops, we end up with the following picture:
\begin{equation} \label{flopdiag}
  \xymatrix{Y \ar@{-->}[rr]^{\chi} \ar[dr]^{\tm} && Y^+ \ar[dl]_{\tm^+}\\
            & X &}
 \end{equation}
(\cite{KoMo}, Theorem~6.14 and Corollary~6.19). Here $Y^+$ is again a terminal Gorenstein threefold with ${-}K_{Y^+}$ big and nef, having $X$ as an anticanonical model. The map $\chi$ is rational and an isomorphism in codimension one. We superscribe any strict transform under $\chi$ with a ``$+$''--sign. Since $K_{Y^+}  + \eps Z^+$ is $\tm^+$--nef, $Z^+$ is $\tm^+$--nef. As above we conclude that $Z^+$ is contractible in $Y^+$.

\vspace{0.2cm}

Lemma~\ref{V} holds for $Y^+$ instead of $Y$ as long as $|F^+|$ is still
spanned on $Y^+$. This need not be the case. Recall that we have chosen $Y$ as a terminal modification of some
$\KQ$--factorialization $X'$ of $X$; in the above
program we might flop some horizontal curves in $Z$, thereby producing a
base locus.

\begin{lemma} \label{flop}
 $|F^+|$ is spanned unless $m = 3$ and $(Z^+, \sO_{Z^+}(Z^+)) = (\PN_2, \sO_{\PN_2}({-}2))$.
\end{lemma}

Here $|F^+|$ restricted to $Z^+$ corresponds to lines through
a given point.

\begin{proof}
Assume that $|F^+|$ is not spanned. Let
 \[\phi^+\colon Y^+ \lra V^+\]
be the divisorial contraction as in (\ref{phi}), contracting $Z^+$. In
order to decide what $Z^+$ is, we
again use the classification from \cite{Cu}. If $Z^+$ maps to a curve and $\fib$ denotes the general fiber, then $Z^+\cdot\fib = -1$ and ${-}K_{Y^+}\cdot\fib = 1$. On $Y^+$ we have
 \begin{equation} \label{zerl}
   {-}K_{Y^+} = Z^+ + B^+ + mF^+.
 \end{equation}
Since $\Bs|F^+| \cap Z^+ \not= \emptyset$ we must have  $F^+\cdot\fib > 0$. From $B^+\cdot\fib \ge 0$ we conclude $0 < m\fib\cdot F^+ \le 2$ which is impossible since $m > 2$.

If $Z^+$ goes to a point, then $(Z^+, \sO_{Z^+}(Z^+))$ is either
$(\PN_2, \sO({-}1))$ or $(\PN_2, \sO({-}2))$ or $(Q_2 \subset \PN_3,
\sO({-}1))$. Near $\Gamma$ the two surfaces $Z$ and $Z^+$ are
isomorphic. With the original pencil on $Z$ we conclude that $Z^+$
contains a smooth rational curve that meets another irreducible curve
in a single point. From $Z^+\cdot\Gamma = -2$ we infer $(Z^+, \sO_{Z^+}(Z^+)) = (\PN_2, \sO({-}2))$.
Then $|F^+|$ restricted to $Z^+$ is a family of lines. Using
${-}K_{Y^+} \cdot \Gamma = m-2$ and the adjunction formula, we find $m
= 3$. The proof of the Lemma is complete.
\end{proof}

Lemma~\ref{V} also holds in the exceptional case $(Z^+, \sO_{Z^+}(Z^+)) = (\PN_2,
\sO_{\PN_2}({-}2))$ for some terminal modification of $X'$, we only
cannot argue as above. Instead, we proceed as follows.

We first run the relative $(K_Y + \eps Z)$--program with respect
to $Y \to X'$, where $X'$ is the above $\KQ$--factorialization of
$X$. In the end we may assume that $Z$ is at least nef on every
$K_Y$--trivial curve contained in a fiber of the pencil $Z \to
\PN_1$. Omitting some details, we conclude that a single flop of a $K_Y$--trivial section of $Z$ transforms
$Y$ into $Y^+$ in (\ref{flopdiag}) and $Z$ into $Z^+ = \PN_2$ as above. Then
 \[Z \simeq \Sigma_1\] 
and $Z\cdot\fib = -1$ for the general fiber $\fib \simeq \PN_1$. We conclude that $Z$ must be negative on at least one extremal ray in $\overline{NE}(Y/\PN_1)$ and conclude
Lemma~\ref{V} as above.

\vspace{0.2cm}

For the proof of Proposition~\ref{cone} it remains to show that $V$ in
Lemma~\ref{V} is a terminal modification of $U_m$ in
section~\ref{examples}. In order to prove this, we consider the system
  \[|{-}K_V + \lambda F_V|, \quad \lambda \ge 0,\]
and choose $\lambda$ such that $m + \lambda \ge 4$. Restricted to $S$ we get $2\Gamma_0 + (m+\lambda)f$, which is now big and nef. Then ${-}K_V + \lambda F_V$ is big and nef and by the Kawamata--Viehweg vanishing theorem $H^1(\sO_V(\lambda F_V)) = H^1(\sO_V(K_V + ({-}K_V + \lambda F_V)) = 0$ implying surjectivity of
  \[H^0(V, \sO_V({-}K_V + \lambda F_V)) \lra H^0(S, \sO_S(2\Gamma_0 + (m + \lambda)f)).\]
Then, since $|F_V|$ is free and $|2\Gamma_0 + (m+\lambda)f|$ is free,
$|{-}K_V + \lambda F_V|$ is free. For $\lambda \ge 1$ and $m + \lambda
\ge 5$, any irreducible curve having zero intersection with ${-}K_V +
\lambda F_V$ must lie in a member of $|F_V|$. This follows immediately
from ${-}K_V + \lambda F_V = ({-}K_V + (\lambda-1)F_V) + F_V =$ nef
$+$ nef. The system is free, for example, if we choose $\lambda = 1$,
for $m \ge 4$, and $\lambda = 2$, for $m = 3$. 

Fix this choice from now on. The map associated to $|{-}K_V + \lambda
F_V|$ is generically 2--to--1 sending $V$ to a variety of minimal degree
  \[\nu\colon V \lra W \subset \PN_{2m + 3\lambda-2}.\]
Since $W$ comes with a pencil $|F_W|$, it must be a scroll. We may rescale
the entries such that ${-}K_V \simeq \nu^*\sO_W(1)$. Then $W \simeq
\PF(d_1, d_2, d_3)$, $d_1 \ge d_2 \ge d_3 \ge -1$, where $d_3 = -1$ in the case $m = 3$, while $d_3 \ge 0$ for $m
\ge 4$. Stein factorization of $V \to W$ leads to a canonical Gorenstein threefold $U$ and a double cover
 \[\mu\colon U \lra W \simeq \PF(d_1, d_2, d_3),\]
such that ${-}K_U = \mu^*\sO_W(1)$. Hence $\mu$ is ramified along a reduced divisor
 \[D \in |\sO_W(4) - 2(d_1 + d_2 + d_3 - 2)F_W|.\]
From $(\sO_W(1))^3 = \frac{1}{2}({-}K_V)^3 = 2m-4$ we infer
  \[d_1 + d_2 + d_3 = 2m-4.\]
The only section of $H^0(V, -K_V-mF_V)$ is the one corresponding to
the image of $B$ in $V$ (cf. (\ref{zbf})). Since
$\mu$ is fiberwise ramified along a quartic, we also have $h^0(W,
\sO_W(1) - mF_W) = 1$, implying
  \[d_1 = m, \quad d_2 < m.\]
In the special case $m = 3$ we have $d_3 = -1$ and $W \simeq \PF(3, 0,
-1)$. It remains to consider the case $m \ge 4$. 

Denote the image of $B$ in $W$ by $B_W$. If $d_3 > 0$, then $2B_W$ is a component of $D$. But $D$ is
reduced, hence we must have $d_3 = 0$. Then $d_1 = m$, $d_2 = m-4$, i.e.,
 \[V \lra U \lra W \simeq \PF(m,m-4,0).\]
We have seen in section~\ref{examples}, that $U = U_m$ can never have canonical singularities for $m \ge 13$, hence $m \le 12$.

Back on the surface $S \in |{-}K_V|$ in Lemma~\ref{V}, we see that
$S$ is generically a double cover of some member $H \in
|\sO_W(1)|$. The map $\nu$ sends $S$ to $\PF(m, m-4)$ and $\Gamma_0$
lies over the minimal section, which is the restriction of the above divisor
$B_W$. In particular, $\Gamma_0$ is not contracted by $V \to U_m$ and
does not meet any curve contracted, i.e., $\Gamma_0 \subset
U_{m,\mathrm{reg}}$ and $V$ is isomorphic to $U_m$ near
$\Gamma_0$. This completes the proof of Proposition~\ref{cone}.
\end{proof}

%%%%%%%%

\section{The Case $W$ a Ruled Surface}
\setcounter{lemma}{0}

This case is as in \cite{Isk}. Instead of $Y$ and $Y_{\Gamma}$ we focus on $X$ and $X_{\Gamma}$, and
diagram (\ref{xgamma}). We use the notation introduced in section~\ref{gensetting}.

\begin{proposition}
In the case $W \simeq \Sigma_{a-b}$, $a > b$, $X$ is the blowup of a sextic in $\PN(1^3, 2,3)$ along an irreducible curve of arithmetic genus one (and $a = 0$, $b = -1$, $m = 3$).
\end{proposition}

\begin{proof}
Since ${-}K_{X_{\Gamma}}$ is ample on $E_{\Gamma}$, we have $b \ge -1$ and $a \ge 0$. Hence
 \[Z_{\Gamma,X}\cdot \ms_{a-b} = b-a < 0 \quad \mbox{and} \quad
 {-}K_{X_{\Gamma}}\cdot \ms_{a-b} = b+2 > 0,\]
where $\ms_{a-b} = E_{\Gamma} \cap Z_{\Gamma,X}$ is the minimal section
of $E_{\Gamma}$. Since $Z_{\Gamma,X}$ is trivial on any
$K_{X_{\Gamma}}$--trivial curve, we conclude that $Z_{\Gamma,X}$ must
be negative on at least one extremal ray in $\{K_{X_{\Gamma}} < 0\}$. Denote by
  \[\phi_X\colon X_{\Gamma} \lra V_X\]
the contraction of this ray. It is a birational map with exceptional
set $Z_{\Gamma,X}$ by \cite{Ben}. Since $Z_{\Gamma,X}$ contains
$K_{X_{\Gamma}}$--trivial curves, it is contracted to a curve.

If $Z_{\Gamma,X}$ is singular along a curve, then its normalization is
a smooth ruled surface. The second map implies that it is $\PN_1
\times \PN_1$. Since $\ms_{a-b} \subset Z_{\Gamma,X,\mathrm{reg}}$ does not
meet the singular locus, we must have 
  \[\deg N_{\ms_{a-b}/Z_{\Gamma,X}} = a = 0,\]
implying $b = -1$. If $Z_{\Gamma,X}$ is smooth in codimension one,
then $h^1(Z_{\Gamma,X}, \sO_{Z_{\Gamma,X}}) \le 1$ by \cite{ReidKaw} and
Iskovskikh's original argument applies: using the ideal sequence of
$Z_{\Gamma,X}$ and the identity $-K_{X_{\Gamma}} = Z_{\Gamma,X} +
(a+2)F_{\Gamma}$ (cf. section~\ref{gensetting}), we see
 \[h^1(Z_{\Gamma,X}, \sO_{Z_{\Gamma,X}}) = h^2(X_{\Gamma}, \sO_{X_{\Gamma}}(-Z_{\Gamma,X})) = h^1(X_{\Gamma}, \sO_{X_{\Gamma}}(-(a+2)F_{\Gamma})).\]
Then the ideal sequence of $(a+2)$ general members of $|F_{\Gamma}|$
 \[0 \lra \sO_{X_{\Gamma}}(-(a+2)F_{\Gamma}) \lra \sO_{X_{\Gamma}}
 \lra \sO_{(a+2)F_{\Gamma}} \lra 0\]
yields $h^0(\sO_{(a+2)F_{\Gamma}}) -1 \le 1$, hence $a \le 0$.

Since $E_{\Gamma}\cdot\ms_{a-b} = a = 0$, the image $Z_X$ of
$Z_{\Gamma,X}$ is still contractible. We can even explicitly give the
supporting divisor: denote the image of $F_{\Gamma}$ in $X$ by $F$. They are Cartier, since $\Gamma \subset X_{\mathrm{reg}}$. The supporting divisor is
 \[H = Z_X + F \in \Pic(X),\]
which is big and nef. Indeed, $\sigma^*H = Z_{\Gamma,X}+F_{\Gamma}+E_{\Gamma}$. Since $Z_{\Gamma,X}+F_{\Gamma} = \varphi^*(\ms_1+\fib)$ is nef, and $\sigma^*H$ restricted to $E_{\Gamma}$ is trivial, $H$ is nef. A direct computation shows $H^3 = 1$. By the Base Point Free Theorem, $|kH|$ is free for $k \gg 0$, defining a birational contraction
 \[\phi\colon X \lra V,\]
contracting $Z_X$ to a curve. The base locus $\Gamma \subset Z_X$ is contracted to a point, the general fiber of the elliptic pencil on $Z_X$ is a section. The variety $V$ is again a Gorenstein Fano threefold with canonical singularities and
 \[K_X = \phi^*K_V+Z_X.\]
From $\phi^*K_V = K_X-Z_X = -2H$ we conclude that ${-}K_V$ is divisible by $2$ in $\Pic(V)$. From $H^0(X, kH) = 1 + \frac{k}{6}(8+3k+k^2)$ we see that $V$ is a sextic in $\PN(1^3,2,3)$.
\end{proof}

\

\begin{proposition}
If $W \simeq \PN_1 \times \PN_1$, then $X \simeq \PN_1 \times S_1$, where $S_1$ denotes a normal del Pezzo surface of degree $1$ (and $a = b = 0$, $m = 4$). 
\end{proposition}

\begin{proof}
In this case, $Z_{\Gamma,X}$ is the pullback of one ruling of $W =
\PN_1 \times \PN_1$. The general fiber of $Z_{\Gamma,X}$ is a smooth
elliptic curve, and $Z_{\Gamma,X}$ meets the singular locus of
$X_{\Gamma}$ at most in points. Going from $X_{\Gamma}$ to $Y_{\Gamma}$, we see
  \[a \le 0.\]
Since $E_{\Gamma} \simeq W$, we have $a = b$, and $X$ Fano implies $a = b = 0$. Since $\varphi$ followed by the natural projection $W \to \PN_1$ contracts all the fibers of $\sigma\colon X_{\Gamma} \to X$ to points, we obtain an induced map
 \[X \lra \PN_1\]
with general fiber $F = \sigma(F_{\Gamma})$ and section $\Gamma$, where $F$ is a normal del Pezzo surface of degree one. We have ${-}K_{X_{\Gamma}} = Z_{\Gamma,X} + 2F_{\Gamma}$. As above,
 \[{-}K_X = Z_X + 2F,\]
and we see that $Z_X$ is nef, so $|kZ_X|$ is free for $k \gg 0$. The map defined by $|kZ_X|$ is a $\PN_1$--bundle with section $F$ and fiber $\Gamma$. As in \cite{Isk} we conclude that $X \simeq F \times \PN_1$ is a product.
\end{proof}

%%%%%%%%%%%%%%%%%%%%%%%

\end{document}